\documentclass[11pt]{amsart}
\usepackage{amsfonts,a4wide}
\usepackage{amssymb}
\usepackage{graphicx}
\usepackage{mathtools}
\usepackage{enumerate}

\newtheorem{theorem}{Theorem}[section]
\newtheorem{lemma}[theorem]{Lemma}
\newtheorem{remark}[theorem]{Remark}
\newtheorem{definition}[theorem]{Definition}
\newtheorem{corollary}[theorem]{Corollary}
\newtheorem{proposition}[theorem]{Proposition}
\newtheorem{example}[theorem]{Example}
\numberwithin{equation}{section}
\numberwithin{table}{section}
\numberwithin{figure}{section}

\newcommand{\im}{\mbox{\rm im\,}}
\newcommand{\re}{\mbox{\rm re\,}}

\newcommand{\diag}{\mbox{\rm diag}}

\newcommand{\Arg}{\mbox{\rm Arg }}

\newcommand{\Comp}{\mathbb{C}}
\newcommand{\eps}{\varepsilon}

\DeclareMathOperator{\ii}{i}

\DeclareMathOperator{\Span}{span}

\newcommand{\eq} [1] {\begin{equation}\label{#1}}
\newcommand{\en} {\end{equation}}
\newcommand {\eqn}  {\begin{eqnarray}}
\newcommand {\enn}  {\end{eqnarray}}
\newcommand {\bstar}    {\begin{eqnarray*}}
\newcommand {\estar}    {\end{eqnarray*}}
\newcommand {\mat}  [1] {\left[\begin{array}{#1}}
\newcommand {\rix}      {\end{array}\right]}
\newcommand{\norm}[1]{\left\| #1 \right\|}
\newcommand{\set}[1]{\left\{ #1 \right\}}

\catcode`@=11     
\font\tenex=cmex10 
\newdimen\p@renwd
\setbox0=\hbox{\tenex B} \p@renwd=\wd0 
\def\bmat#1{\begingroup \m@th
  \setbox\z@\vbox{\def\cr{\crcr\noalign{\kern2\p@\global\let\cr\endline}}%
    \ialign{$##$\hfil\kern2\p@\kern\p@renwd&\thinspace\hfil$##$\hfil
      &&\quad\hfil$##$\hfil\crcr
      \omit\strut\hfil\crcr\noalign{\kern-\baselineskip}%
      #1\crcr\omit\strut\cr}}%
  \setbox\tw@\vbox{\unvcopy\z@\global\setbox\@ne\lastbox}%
  \setbox\tw@\hbox{\unhbox\@ne\unskip\global\setbox\@ne\lastbox}%
  \setbox\tw@\hbox{$\kern\wd\@ne\kern-\p@renwd\left[\kern-\wd\@ne
    \global\setbox\@ne\vbox{\box\@ne\kern2\p@}%
    \vcenter{\kern-\ht\@ne\unvbox\z@\kern-\baselineskip}\,\right]$}%
  \null\;\vbox{\kern\ht\@ne\box\tw@}\endgroup}
\catcode`@=12 

\def\Span{\mathop{\mathrm{span}}}
\def\diag{\mathop{\mathrm{diag}}}

\newcommand {\comment}[1]{} 
\usepackage{color}

\title{Stability of Matrix Polynomials in One and Several Variables}
\author{Oskar Jakub Szymański, Micha\l{} Wojtylak}
\date{\today}

\begin{document}

\maketitle
\begin{abstract}
The paper  presents methods of eigenvalue localisation of regular matrix polynomials, in particular, stability of matrix polynomials is investigated. For this aim a stronger notion of hyperstability  is introduced and widely discussed.
Matrix versions of the Gauss-Lucas theorem and Sz\'asz inequality are shown.  Further,  tools for investigating (hyper)stability   by multivariate complex analysis methods are provided. Several second- and third-order matrix polynomials with particular semi-definiteness  assumptions on coefficients are shown to be stable.  
\end{abstract}

\begin{footnotesize}
    \begin{align*}
 &\textbf{Keywords} : \text{matrix polynomial, eigenvalue, stability, polarization operator, multivariate polynomial.} \\
 &\textbf{MSC2020} : \text{15A18, 15A22, 15B57, 30C15, 32A60}
    \end{align*}
\end{footnotesize}

\section{Introduction}


Given a matrix polynomial $P(\lambda)=\lambda^d A_d+ \lambda^{d-1} A_{d-1} + \dots + A_0$ ($A_0, A_1, \dots , A_d \in \Comp^{n, n}$) it is a natural question  to ask under which conditions all eigenvalues are located outside a given set $D$. By eigenvalues we mean zeros of the function $\lambda\mapsto\det P(\lambda)$, which is  assumed to be nonzero, i.e.,  $ P(\lambda)$ is regular. Note that $ P(\lambda)$ has no eigenvalues in some set $D$ if and only if  
\begin{equation}\label{SC}
  \textrm{ for all } {\mu \in D}, \;{\textrm{ for all } } {x \in \mathbb{C}^n \setminus \{0\}},\; \textrm{ there exists } {y \in \mathbb{C}^n } \textrm{ such that } y^*P(\mu)x \neq 0.
\end{equation}
 There are several techniques proving such localisation of eigenvalues. One of them is the numerical range \cite{LiR94, Psa03, Psa00}, see Definition~\ref{numran} below, (see also \cite{bini2013, hig2003} for other techniques less related with the current project). Indeed, it is elementary  that the numerical range contains all eigenvalues, cf.  \cite{mar1997, MehMW22} for applications.  However, a localisation condition for eigenvalues that uses numerical range is in general very restrictive.
Consider for example a port-Hamiltonian pencil of the form (cf.\cite{MehMW18}, see also \cite{GerHH21,mehrmann2022control})
\begin{equation}\label{pH}
\lambda E+(J+R)Q
\end{equation}
with matrices $J$, $E^*Q$ positive semi-definite, and $J$ being skew-symmetric and $Q$ not necessarily invertible, but satisfying a weaker condition. It was shown in \cite{MehMW18} that the eigenvalues are in the closed left half-plane.  The key issue of that proof, although not stated explicitly,  was considering the numerical range of  the pencil 
$\lambda Q^* E+Q^*(J+R)Q,
$ not of the original one.  
 In order to cover such situations, we introduce the following condition
 \begin{equation}\label{HSC}
 \textrm{ for all }  x \in \mathbb{C}^n \setminus \{0\},\;\textrm{ there exists } y \in \mathbb{C}^n, \textrm{ such that for all } \mu \in D : y^*P(\mu)x \neq 0.
 \end{equation}
  We will refer to condition \eqref{HSC} as to \textit{hyperstability  with respect to $D$}, see Definition \ref{defla}. Note that on one hand it is clearly stronger than \eqref{SC}, on the other it is clearly weaker than ordering that the numerical range is outside $D$ (take $y=x$).  The choice $y=Qx$ shows that the pencil in \eqref{pH} is hyperstable with respect to the open right half-plane, this idea will be extended in the current paper to quadratic polynomials in Corollary~\ref{c-hp}.

 The first outcome of the present manuscript is presenting some classes of matrix polynomials for which the conditions \eqref{SC} and \eqref{HSC} are equivalent, see Theorem \ref{uppert}, and providing examples showing that they are not equivalent in general, cf. Examples \ref{exa} and \ref{hyper-nsinf}. In particular, for linear pencils the two notions coincide. Further, we show that similarity of matrix polynomials do not preserve hyperstability (Remark~\ref{eqstab}).  For these two reasons, we think that the notion does not seem to be much related with invariant factors, and we concentrate the current research on location of eigenvalues only. On the other hand, the notion is related with current work on (generalised) triangularisation of matrix polynomials  \cite{ang2021, trian2013, tis2013}, see Proposition \ref{upperblock} and Theorem \ref{uppert}.
 
Further, it appears that the hyperstability notion in \eqref{HSC} inherits  many properties of  stability for scalar polynomials. For example, we show a hyperstable version of the Gauss-Lucass Theorem for matrix polynomials (Theorem~ \ref{GLmat}), while  an analogous result using the notion of stability \eqref{SC} simply fails, see Examples \ref{nonGL}.   

Our next contribution is introducing several complex variables. Up to our knowledge, the notions of regularity and  numerical range for multivariate matrix polynomials are  not yet established and need a separate investigation. We show that the condition of hyperstability  may be used as an alternative,  as
the condition \eqref{HSC}  generalises naturally onto multivariate matrix polynomials, see Definition~\ref{defzi}. Further, we introduce the polarization operator for matrix polynomials, see Section~\ref{sPol}. E.g., for a quadratic matrix polynomial $\lambda^2 A_2 +\lambda A_1+ A_0$ its polarization is given by $z_1z_2A_2+\frac{z_1+z_2}2A_1+A_0$.
It is well known that the polarization operator preserves stability of scalar polynomials of any degree.
This concept was heavily used e.g. in \cite{BorB09}, where all operators preserving stability of multivariate polynomials were classified. 
It appears that polarisation preserves also by hyperstability of matrix polynomials, see Theorem~\ref{GWS}. Again, as in the case of the Gauss-Lucas Theorem,  an analogue result, with stability condition instead of hyperstability, is not true (see Example \ref{nonstab}).

  Another outcome of the paper are general hyperstability criteria for quadratic and cubic matrix polynomials, see Theorems \ref{poly2} and \ref{poly3}. For example, we show that if the function 
  $\det(z_1^2 A_2+z_2 A_1 +A_0)$ has no zeros in some  set $D^2$ ($D\subseteq\Comp$) then the matrix polynomial
  $\lambda^2 A_2+\lambda A_1+A_0$ satisfies the hyperstability condition \eqref{HSC} with respect to $D$. 
 Note that analytic properties of a certain multivariate polynomial imply here strictly linear algebra properties of a matrix polynomial $P(\lambda)$, which we find remarkable.

The last outcome of the paper is proving hyperstability and stability of some one variable matrix polynomials of degrees 2 and 3 and certain multivariate matrix polynomials. The methods are particularly tailored for matrix polynomials with certain positivity conditions on the coefficients. Such polynomials were recently considered in \cite{MehMW21,MehMW22} due their numerous applications in mathematical modelling \cite{BetHMST08,MehMW21,mehrmann2022control}. In particular, Theorems \ref{half-plane}, \ref{ker} and \ref{deg3} contribute to the knowledge of the topic.

The paper is organised as follows.  Section~\ref{sPrel} contains the usual linear algebra notions and notations.  In Section~\ref{sHyper} we introduce  hyperstability and show its basic properties.  In Section~\ref{sGL} we prove a generalisation of the Gauss-Lucas  theorem for hyperstable matrix polynomials and a version of a Sz\'asz-type inequality for matrix polynomials. In Section~\ref{s5} we introduce the complex analysis of many variables and show some hyperstability criteria for cubic and quadratic matrix polynomials.
In Section~\ref{sPol} we introduce the polarisation operator and show that it preserves hyperstability. 
The last part of the paper, Section~\ref{sPos}, states hyperstabilty for particular classes of quadratic and cubic matrix polynomials  with coefficients satisfying certain positive-definiteness conditions.



\section{Preliminaries}\label{sPrel}

Let us begin with fixing the notation. 
The matrix polynomials of one variable will be denoted as
$P(\lambda) = A_d\lambda^d + A_{d-1}\lambda^{d-1} + \dots +A_0$,  with matrix coefficients $A_0, A_1, \dots , A_d\in\Comp^{n,n}$ and $\lambda$ being the complex variable. By multivariate matrix  polynomials we will understand finite sums of the form $P(z_1, z_2, \dots, z_\kappa)= \sum_{\alpha_1, \alpha_2, \dots , \alpha_\kappa} z_1^{\alpha_1}z_2^{\alpha_2} \cdots z_\kappa^{\alpha_\kappa} A_{\alpha_1, \alpha_2, \dots , \alpha_k}$, where $A_{\alpha_1, \alpha_2, \dots , \alpha_k} \in \Comp^{n,n}$, $\alpha_1, \alpha_2 \dots , \alpha_\kappa \geq 0$, and   $z_1, z_2 \dots , z_\kappa$ are complex variables.
Although we allow $\kappa=1$, this will usually lead to a trivial situation.

For scalar polynomials ($n=1$) we will use lower-case letters for coefficients. Let $D$ be a nonempty open or closed subset of $\Comp^\kappa$. We will call a scalar polynomial   $p \in \mathbb{C}[z_1, z_2, \dots , z_{\kappa}]$ \emph{stable} with respect to $D$  if and only if it has no zeros in $D$.
Two possible generalisations of this notion to matrix polynomials, see Definitions~\ref{defla} and \ref{defzi}, are the main object of this paper.

    \begin{definition}
\rm We call a matrix polynomial $P(\lambda) = A_d\lambda^d + A_{d-1}\lambda^{d-1} + \dots +A_0 \in \mathbb{C}^{n, n}[\lambda]$, where $A_d \neq 0$, \emph{regular} if and only if the function $\lambda \mapsto \det P(\lambda)$ is a nonzero scalar polynomial. In such case, we call a vector $x \in \mathbb{C}^n \setminus \{0\}$ \emph{eigenvector} of the matrix polynomial $P(\lambda)$ if and only if there exists $\lambda_0 \in \mathbb{C}$ such that
		\begin{equation*}
P(\lambda_0)x = 0\text{.}
		\end{equation*}
Then, such complex number $\lambda_0$ is called an \emph{eigenvalue} corresponding to the eigenvector $x$. Note that for a regular polynomial being an eigenvalue is equivalent to being a root of the scalar polynomial $\lambda \mapsto \det P(\lambda)$.
    \end{definition}
Note that we do not intend to define regular or singular multivariate polynomials,
the example $z_1(A+z_2I)$ shows difficulties with generalising the notion (fixing $z_1=0$ gives a singular polynomial in the second variable). Instead, we will heavily use the hyperstability notion,  see Definition \ref{defzi}. Multivariate matrix polynomials will appear again in Section~\ref{s5}, now we concentrate on scalar ones.
In few cases, we need a factorization tools like the well known Kronecker canonical form of a matrix pencil or Smith canonical form of an arbitrary square matrix polynomial reviewed below, see \cite[Chapter VI]{Gan59}.
    \begin{definition}
\rm Let $P(\lambda) \in \mathbb{C}^{n, n}[\lambda]$ be a regular matrix polynomial. Then, there exist regular polynomials $U(\lambda), V(\lambda) \in \mathbb{C}^{n, n}[\lambda]$ with nonzero constant determinants such that
        \begin{equation}\label{Smith}
U(\lambda)P(\lambda)V(\lambda) = \diag{\Bigl(s_1(\lambda), s_2(\lambda), \dots , s_n(\lambda)\Bigr)} =: S(\lambda)\text{,}
        \end{equation}
where the symbols $s_1(\lambda), s_2(\lambda), \dots , s_n(\lambda) \in \mathbb{C}[\lambda]$ denote uniquely determined monic polynomials such that a polynomial $s_i(\lambda)$ is a divisor of a polynomial $s_{i+1}(\lambda)$ for $i \in \{1, 2, \dots , n\}$. The non-singular diagonal matrix polynomial $ S(\lambda) \in \mathbb{C}^{n, n}$ is called \emph{the Smith canonical form} of  $ P(\lambda)$ and the diagonal entries in the right hand side of \eqref{Smith} are called \emph{the invariant factors} of $ P(\lambda)$.
    \end{definition}

Besides the notion of an eigenvalue of a matrix polynomial $P \in \mathbb{C}^{n, n}[\lambda]$ we need a notion of a wider set containing all eigenvalues, i.e. the notion of the numerical range. It was introduced in \cite{LiR94} and studied further in \cite{Psa03,Psa00}.
    \begin{definition}\label{numran}
\rm Let $P(\lambda) \in \mathbb{C}^{n, n}[\lambda]$ be an arbitrary matrix polynomial. A subset
        \begin{equation}
W(P) := \{\mu \in \mathbb{C} : x^*P(\mu)x = 0 \text{ for some } x \neq 0\}
        \end{equation}
of the complex plane $\mathbb{C}$ is called \emph{the numerical range} of the matrix polynomial $P(\lambda)$.
    \end{definition}

\section{Hyperstability of matrix polynomials in one variable}\label{sHyper}

In this section we will deal with matrix polynomials of one variable, 
we emphasize now two central notions.
	
	\begin{definition}\label{defla}
\rm Let the symbol $D$ denote an nonempty open or closed subset of the complex plane $\mathbb{C}$ and let $P(\lambda) \in \mathbb{C}^{n, n}[\lambda]$ be a matrix polynomial. We say that the polynomial $ P(\lambda)$ is \emph{stable with respect to $D$} if and only if it does not have eigenvalues in $D$.
Further, we say that the polynomial $P(\lambda)$ \emph{is hyperstable with respect to $D$} if and only if for all $x \in \mathbb{C}^n \setminus\{0\}$ there exists $y \in \mathbb{C}^n \setminus \{0\}$ such that 
		\begin{equation}\label{nozeroinD}
y^*P(\mu)x \neq 0 \;\;\text{for all} \;\;\mu \in D\text{.}
		\end{equation}
	\end{definition}
		
Observe that stability of a matrix polynomial implies its regularity. Further, the notion of hyperstability is situated somewhere between stability and a (rather strong) numerical range condition.
\begin{proposition}\label{abc}
Let $P(\lambda) \in \mathbb{C}^{n, n}[\lambda]$ be a matrix polynomial and let $D$ be a nonempty open or closed subset of the complex plane $\mathbb{C}$. Consider the following conditions.
\begin{enumerate}[\rm (a)]
\item\label{num}  the numerical range $W(P)$ does not intersect $D$;
\item\label{hyper} the polynomial $P(\lambda)$ is hyperstable with respect to $D$; 
\item\label{stable} the polynomial $ P(\lambda)$ is stable with respect to $D$.
\end{enumerate}
Then \eqref{num}$\Rightarrow$\eqref{hyper}$\Rightarrow$\eqref{stable}.
\end{proposition}
		
\begin{proof}
The first implication follows by setting $y=x$ in \eqref{nozeroinD}. The second implication becomes obvious, as one observes that stability can be reformulated as:
for all $\mu\in D$ and for all $x \in \mathbb{C}^n \setminus\{0\}$ there exists $y \in \mathbb{C}^n \setminus \{0\}$ such that $y^*P(\mu)x \neq 0$.
\end{proof}		

 It appears that in many cases implication \eqref{num} $\Rightarrow$ \eqref{stable} above is a convenient criterion for stability, see \cite{MehMW22}. Note that the numerical range of a polynomial $P(\lambda)=\lambda I_n - A$ coincides with the numerical range of the matrix $A$. Thus condition \eqref{num} cannot be, in general, equivalent to stability \eqref{stable}. As for matrix pencils hyperstability \eqref{hyper} is equivalent to 
stability \eqref{stable} (see Theorem~\ref{uppert}\eqref{pq} below), we obtain that in many cases \eqref{num} is not  equivalent to \eqref{hyper}. Further, hyperstability \eqref{hyper} is,  in general, not equivalent to stability \eqref{stable}, cf. the important following Example.

\begin{example}\label{exa}\rm	
Let
$$
P(\lambda) := \mat{cc} 1 & \lambda \\ \lambda & \lambda^2 + 1\rix.
$$	
Then $\det P(\lambda)\equiv 1$, hence $P(\lambda)$ is a stable matrix polynomial with respect to any open or closed subset $D$ of the complex plane.
However, taking $x=[0 \; 1]^{\top}$ and arbitrary $y = [y_1 \; y_2]^{\top}$ yields 
$$
y^*P(\lambda)x = \overline{y}_1 \lambda + \overline{y}_2 (\lambda^2 + 1)\text{.}  
$$	
Note that for any $y \neq 0$ the polynomial above has a root in the closed unit disc $\overline{\mathbb{D}}$, as either $y_2 = 0$ and we have one root $\lambda=0$ or $y_2\neq 0$ and we have two roots (counting multiplicities) with the product equal $1$.  Replacing $\lambda$ by $\alpha\lambda + \beta$ one gets an example of a polynomial stable, but not hyperstable, with respect to an arbitrary desired half-plane or disc. E.g., the polynomial $P(\lambda-\ii)$  is stable, but not hyperstable, with respect to the half-plane $H_0$.

Observe, that the polynomial $P(\alpha\lambda+\beta)$ ($\alpha\neq0$) has always infinity as an eigenvalue, i.e., the leading coefficient $\diag(0,\alpha^2)$ is not an invertible matrix. One may wonder if this property is necessary to construct such an example. Apparently, this is not the case, see Example~\ref{hyper-nsinf} below.
\end{example}

Let us recall two notions of equivalence for matrix polynomials, cf. \cite{Gan59}.
We say that matrix polynomials $P(\lambda), Q(\lambda) \in \mathbb{C}^{n, n}[\lambda]$ are \emph{equivalent} if and only if there exist  matrix polynomials  $U(\lambda), V(\lambda) \in \mathbb{C}^{n, n}[\lambda]$ with constant nonzero determinants such that $P(\lambda) = U(\lambda)Q(\lambda)V(\lambda)$. We say that the polynomials $ P(\lambda)$ and $Q(\lambda)$ are \emph{strongly equivalent} if and only if the polynomials $ U(\lambda)$ and $ V(\lambda)$ are constant invertible matrices. Clearly, the latter relation preserves hypertability, in fact a stronger statement holds: 
\begin{lemma}\label{lQ}
Let $D$ be any open or closed set, let $Q\in\Comp^{n,n}$ be any matrix and let $S\in\Comp^{n,n}$ be invertible. If $Q^* P(\lambda)S$  is hyperstable with respect to $D$,  then $P(\lambda)$ is hyperstable with respect to $D$. 
\end{lemma}

 The proof follows directly from the definition. In a moment we will see that  equivalence does not preserve hyperstability, see Remark~\ref{eqstab}. First, however, we need to show some properties of hyperstability, connected with (block) upper-triangular matrices. 

\begin{proposition}\label{upperblock}
Let $P(\lambda) \in \mathbb{C}^{n,n}[\lambda]$ be a matrix polynomial and let $D$ be an open or closed subset of the complex plane. Assume that  $ P(\lambda)$ is strongly equivalent to a block upper-triangular matrix polynomial 
$$
 \mat{ccc} P_{11}(\lambda) & \cdots & P_{1m}(\lambda)\\ &\ddots&\vdots\\ && P_{mm}(\lambda)  \rix,\quad P_{ij}(\lambda)\in\Comp^{k_i,k_j}[
 \lambda], \ \sum_{j=1}^m k_j=n,
$$
with the diagonal entries hyperstable with respect to $D$. Then the polynomial $P(\lambda)$  is hyperstable with respect to $D$.
\end{proposition}
	
	\begin{proof}
	As strong equivalence preserves hyperstability, we may assume without loss of generality that $U =V = I_n$. Take $x = \mat{c} x_1^\top  \dots x_m^\top \rix^\top \neq 0$,  with 
	$x_j\in \Comp^{k_j}$ and let $r \in \set{1, 2, \dots , n}$ denote the index of the last nonzero  $x_j$. Let $y_r\in\Comp^{k_r}\setminus\set0$ be such that the polynomial $ y_r^* P_{rr}(\lambda)x_r$ is stable with respect to $D$. Taking $y=\mat{c} 0 \dots 0\ y_r^\top 0 \dots 0 \rix^\top  $  we get  $y^* P(\lambda)x=y_r^* P_{rr}(\lambda)x_r$, which  is stable with respect to $D$.  
	
	\end{proof}
	
	The following Theorem presents a  class of matrix polynomials for which the notions of stability and hyperstability coincide. For applications  see Proposition~\ref{MGT} below.

\begin{theorem}\label{uppert}
Let $P(\lambda) \in \mathbb{C}^{n, n}[\lambda]$ be a matrix polynomial and let $D$ be an open or closed subset of the complex plane, then the following holds.
\begin{enumerate}[\rm (i)]
\item\label{ut} Assume that a matrix polynomial $P(\lambda)$ is strongly equivalent to an upper triangular matrix polynomial. Then it is hyperstable with respect to $D$ if and only if it is stable with respect to $D$.
\item\label{pq} If $P(\lambda) = p(\lambda)A + q(\lambda)B$ with some scalar polynomials $p(\lambda), q(\lambda) \in \mathbb{C}[\lambda]$ and $A, B \in \Comp^{n, n}$ (in particular: if $P(\lambda)$ is a matrix pencil, i. e. a matrix polynomial of degree one), then it is hyperstable with respect to $D$ if and only if it is stable with respect to $D$.
\end{enumerate}
\end{theorem}
	
	\begin{proof}
\eqref{ut} Consider a matrix polynomial $P(\lambda)$ stable with respect to $D$. We have 
$$
U^{-1}P(\lambda)V^{-1} = \mat{ccc} p_{11}(\lambda) & \cdots & p_{1n}(\lambda)\\ &\ddots&\vdots\\ && p_{nn}(\lambda)  \rix
$$
with a stable upper triangular matrix polynomial on the right hand side of the equation. Hence, the scalar polynomials on the diagonal $p_{11}(\lambda), p_{22}(\lambda), \dots , p_{nn}(\lambda) \in \mathbb{C}[\lambda]$ are stable with respect to $D$. By Proposition~\ref{upperblock} $P(\lambda)$ is hyperstable with respect to $D$.   
	
\eqref{pq} Using the Kronecker form \cite{Gan59} (or the generalised Schur form) of matrices $A$ and $B$ we obtain the matrix polynomial  $\lambda \mapsto p(\lambda)A + q(\lambda)B$ to be strongly equivalent to an upper triangular one and the claim follows from \eqref{ut}.
	\end{proof}

\begin{remark}\label{eqstab}\rm
Note that an arbitrary matrix polynomial $P(\lambda) \in \mathbb{C}^{n, n}[\lambda]$ is equivalent to an upper triangular matrix polynomial (with $V(\lambda) = I_n$), see \cite[Chapter VI]{Gan59}. E.g., the polynomial in Example~\ref{exa} satisfies the following equality:
$$
\mat{cc} 1 & \lambda\\ \lambda & \lambda^2 + 1\rix = \mat{cc} 1 & 0\\ \lambda & 1\rix\mat{cc} 1 & \lambda\\ 0 & 1\rix\text{.}
$$
Hence, equivalence of matrix polynomials does not preserve hyperstability.  
\end{remark}	

\section{Gauss-Lucas Theorem and Sz\'asz type inequality for matrix polynomials}\label{sGL}

We present here generalisations of  two major results for scalar polynomials onto matrix polynomials. 
The celebrated Gauss-Lucas Theorem says that for a non-constant scalar polynomial $p(\lambda)$ the roots of 
$p'(\lambda)$ are contained in the convex hull of the roots of $p(\lambda)$. Using the terminology of the current paper one can formulate it as follows: if $D\subseteq\Comp$ is such that $\Comp\setminus D$ is convex and $p(\lambda)$ is stable with respect to $D$ then $p'(\lambda)$ is stable with respect to $D$. It seems to be generally  known, that a similar statement cannot be true for matrix polynomials. Nonetheless,  we present a very general example which will serve in future constructions.

\begin{example}\label{nonGL}\rm
Consider a polynomial 
$$
P(\lambda)=\mat{cc} \lambda p(\lambda)+q(\lambda) & p(\lambda)\\ \lambda & 1\rix,
$$
where $q(\lambda)$ has its roots outside $D$ and $p'(\lambda)$ is non-constant and has some of its roots inside $D$. Then
$$
\det P(\lambda)=q(\lambda),\quad \det P'(\lambda) = -p'(\lambda),
$$
i.e. the polynomial $P(\lambda)$ is stable with respect to $D$ but its derivative $P'(\lambda)$ is not.
\end{example}

The example above says that, in general, there is no relation between the location of eigenvalues of a polynomial $P(\lambda)$ and of its derivative.
However, a hyperstable version of the Gauss-Lucas Theorem holds. 

\begin{theorem}\label{GLmat}
 Let $D\subseteq\Comp$ be such that $\Comp\setminus D$ is convex. If a matrix polynomial $P(\lambda)$ is hyperstable with respect to $ D$ and the entries of its derivative $P'(\lambda)$ are linearly independent polynomials, then the matrix polynomial $P'(\lambda)$ is also hyperstable with respect to $D$.  
\end{theorem}

\begin{proof}
Fix a nonzero $x\in\Comp^n\setminus\{0\}$. As $P(\lambda)$ is hyperstable with respect to $D$ there exists $y\in\Comp^n\setminus\{0\}$ such that the scalar polynomial $p(\lambda)=y^*P(\lambda)x$ has all its roots in $\Comp\setminus D$. Note that the polynomial $p'(\lambda)=y^*P'(\lambda)x$ is a nontrivial linear combination of entries of $P'(\lambda)$, hence $p'(\lambda)$ is a nonzero polynomial.
As $\Comp\setminus D$ is convex, by the Gauss-Lucas theorem, $p'(\lambda)$ has all its roots in $D$.
\end{proof}

Firstly,  note that the assumption on the derivative $P'(\lambda)$ having independent entries cannot be dropped:

\begin{example}\rm
The polynomial $P(\lambda)=\diag(\lambda, 1)$ is hyperstable with respect to the outside of the unit disc $\mathbb{D}$, but its derivative $P'(\lambda)$ is singular, hence is not hyperstable with respect to any nonempty set.  
\end{example}

Secondly, note that the assumption on $P(\lambda)$ being hyperstable cannot be relaxed to stability, even if we keep the assumption on entries of $P'(\lambda)$ being linearly independent and relax the claim to $P'(\lambda)$ being stable with respect to $D$.

\begin{example}\label{hyper-nsinf}\rm
We specify and modify the polynomial $P(\lambda)$ from Example~\ref{nonGL}. Let $D=\{z\in\Comp:|z|>1\}$, $q(\lambda)=\lambda^2$, $p(\lambda)=\lambda^3-4\lambda$ and consider a perturbed polynomial
$$
P_\eps(\lambda):= 
\mat{cc} \lambda p(\lambda)+q(\lambda) & p(\lambda)\\ \lambda & 1+\eps \lambda^4\rix= \mat{cc} \lambda^4-3\lambda^2 & \lambda^3-4\lambda\\ \lambda & 1+\eps \lambda^4\rix.
$$
Fix  $\eps>0$ sufficiently small, so that  $P_\eps(\lambda)$ has still its eigenvalues inside the unit disc $\mathbb{D}$ and $P'_\eps(\lambda)$ has still its eigenvalues outside the closed unit disc. This is possible due to continuity of zeros of polynomials with respect to the coefficients.  Further, note that
$$
P'_\eps(\lambda) =  \mat{cc} 4\lambda^3-6\lambda & 3\lambda^2-4\\ 1 & 4\eps \lambda^3\rix,
$$
which clearly has linearly independent entries. 

By Theorem~\ref{GLmat} the polynomial $P_\eps(\lambda)$ is not hyperstable with respect to $D$.
Further, note that the leading coefficient of this polynomial is an invertible matrix $\diag(1, \eps)$, which gives the desired example of stable, but not hyperstable polynomial with no eigenvalues at infinity. 
\end{example}

We use the current moment to present a Sz\'asz-type inequality for stable matrix polynomials. 
\begin{proposition}\label{Szasz}
Consider a matrix polynomial $P(\lambda) = \lambda^d A_d + \lambda^{d-1} A_{d-1} + \dots + \lambda A_1 + I_n$. 
If the numerical range $W(P)$ is contained in some  half-plane $H_\varphi$, $\varphi \in [0;2\pi)$, then
$$
\| P(\lambda) \| \leq 2\exp \left( \lambda_{H}\Bigl[\lambda A_1 -|\lambda|^2 A_2\Bigr] +\frac12 |\lambda|^2 \| A_1\|^2       \right), \quad \lambda \in \Comp,
$$
where the symbol $\lambda_{H}(X)$ denotes the largest (possibly negative) eigenvalue of the Hermitian matrix $\frac{X+X^*}2$. 
\end{proposition}

\begin{proof}
Without loss of generality we may assume that the numerical range of $P(\lambda)$ is contained in the open upper half-plane. By the well known fact that the norm of a matrix is less than twice its numerical radius, see e.g. \cite{HorJ91},  we obtain for $\lambda\in\Comp$
\begin{eqnarray*}
\|P(\lambda)\| &\leq & 2\sup_{\|x\| = 1} |x^*P(\lambda)x| \\
& \leq & 2\sup_{\|x\|=1} \exp \left(\re(x^*A_1x \lambda) - \re(x^* A_2 x) |\lambda|^2  + \frac12 |x^*A_1x|^2|\lambda|^2  \right),
\end{eqnarray*}
where the second inequality follows from scalar inequality for stable polynomials (\cite[Lemma 5]{deB61}, see also \cite[Theorem 1.3]{Kne19}), applied to $x^*P(\lambda)x$.
Taking the supremum under the exponent we see that the assertion follows.
\end{proof}

In view of Proposition~\ref{abc} one may ask whether hyperstability implies any Sz\'asz-type inequality.

\section{Hyperstability of matrix polynomials in one variable via multivariate matrix polynomials}\label{s5}

In this Section we show how certain analytic properties of multivariate matrix polynomials imply algebraic properties of one variable matrix polynomials. To be more precise, we show that stability of multivariate matrix polynomials (with respect to some set $D^\kappa$) implies hyperstablity of matrix polynomials of one variable. First, let us extend Definition~\ref{defla} to the multivariate case.

\begin{definition}\label{defzi}
\rm Let the symbol $D$ denote a nonempty open or closed subset of the complex plane. We say that a multivariate matrix polynomial $(z_1, z_2, \dots , z_\kappa) \mapsto P(z_1, z_2, \dots , z_\kappa)$ \emph{is stable with respect to $D^{\kappa}$} if and only if for all $(\mu_1,\dots,\mu_\kappa) \in D^{\kappa}$ and all $x \in \mathbb{C}^n\setminus\{0\}$ we have
		\begin{equation*}
P(\mu_1, \mu_2, \dots, \mu_{\kappa})x \neq 0\text{.}
		\end{equation*}
Further, we say that $(z_1, z_2, \dots , z_{\kappa}) \mapsto P(z_1, z_2, \dots , z_{\kappa})$ is \emph{hyperstable with respect to $D^{\kappa}$} if and only if for all $x \in \mathbb{C}^n\setminus\{0\}$ there exists $y \in \mathbb{C}^n\setminus\{0\}$ such that 
		\begin{equation}\label{nozeroinDmatrix}
y^*P(\mu_1, \mu_2, \dots , \mu_{\kappa})x \neq 0 \;\;\text{for all} \;\;(\mu_1, \mu_2, \dots , \mu_{\kappa}) \in D^{\kappa}\text{.}
		\end{equation}
\end{definition}
	
We present now one of the major outcomes of the current paper relating hyperstability of one variable quadratic polynomials to stability of certain two variables matrix polynomials.

\begin{theorem}\label{poly2}
Let $P(\lambda) = \lambda^2 A_2 + \lambda A_1 + A_0$ be a  quadratic matrix polynomial and let $D$ be a nonempty open or closed subset of the complex plane $\mathbb{C}$. If at least one of the following conditions hold:
	\begin{enumerate}[\rm (a)]
\item\label{0?D} the multivariate matrix polynomial  $(z_1, z_2) \mapsto z_1^2 A_2 + z_2 A_1 + A_0$ is stable with respect to $D^2$,
\item\label{0notinD1} the  multivariate matrix polynomial $(z_1, z_2) \mapsto z_1z_2 A_2 + z_2 A_1 + A_0$  is stable with respect to $D^2$ and $0 \notin D$,
\item\label{0notinD2} the  multivariate matrix polynomial  $(z_1, z_2) \mapsto z_1^2z_2 A_2 + z_1^2 A_1 + z_2 A_0$ is stable with respect to $D^2$ and $0 \notin D$,
	\end{enumerate}
then the matrix polynomial $\lambda \mapsto P(\lambda)$ is hyperstable with respect to $D$.
\end{theorem}

\begin{proof}
First observe that setting $z_1 = z_2 = \lambda$ implies, in each case, that the polynomial $\lambda \mapsto P(\lambda)$ is stable with respect to $D$, in particular it is regular. In the case \eqref{0notinD2} stability of a polynomial $\lambda \mapsto \lambda P(\lambda)$ is equivalent to stability of the polynomial $\lambda \mapsto P(\lambda)$. To show that the polynomial $\lambda \mapsto P(\lambda)$ is hyperstable with respect to $D$, fix $x \in \Comp^n\setminus\{0\}$. Below by $\perp$ we denote orthogonality with respect to the standard complex inner-product. The proof in each case is similar, however, requires certain adaptations and thus we present all details.

%

Assume \eqref{0?D} holds. If there exists a vector $y \in \mathbb{C}^n\setminus\{0\}$ such that $y \perp A_2x, \; y \perp A_1x$ and $y \not\perp A_0x$, then 
	\begin{equation*}
y^* P(\lambda)x = y^*A_0x
	\end{equation*}
and \eqref{nozeroinD} is simply satisfied.

Hence, further we assume that such vector $y$ does not exist, i.e.
$
\{A_2x, A_1x\}^{\perp} \subseteq \{A_0x\}^{\perp}$,  equivalently $A_0x\in \Span\set{A_2x, A_1x}$.
Thus, we have 
\begin{equation}\label{A_0x}
A_0x = \alpha_0 A_2x + \beta_0 A_1x
\end{equation}
for some $\alpha_0, \beta_0 \in \Comp$. 
Next, we can write
	\begin{align*}
P(\lambda)x = \lambda^2 A_2x + \lambda A_1x + (\alpha_0 A_2x + \beta_0 A_1x) = (\lambda^2 + \alpha_0)A_2x + (\lambda + \beta_0)A_1x\text{.}
	\end{align*}
	
We show now that at least one of the scalar polynomials $\lambda \mapsto \lambda^2 + \alpha_0$ or $\lambda \mapsto \lambda + \beta_0$ is stable with respect to $D$.	
By the assumption that the polynomial $(z_1, z_2) \mapsto z_1^2 A_2 + z_2 A_1 + A_0$ is stable with respect to $D^2$, we have that the equation $(z_1^2 A_2 + z_2 A_1 + A_0)x = 0$ has no solutions $(z_1, z_2) \in D^2$. Substituting  \eqref{A_0x} we obtain that the equation
$$
(z_1^2 + \alpha_0)A_2x + (z_2 + \beta_0)A_1x = 0
$$
has no solutions $(z_1, z_2) \in D^2$. This implies that at least one of the scalar polynomials 
$\lambda \mapsto \lambda^2 + \alpha_0$ or $\lambda \mapsto \lambda + \beta_0$ has no roots in $D$ and the claim follows.
%

If the polynomial $\lambda \mapsto \lambda^2 + \alpha_0$ is stable with respect to $D$, then similarly as before we seek for a vector $y$ such that $y \not\perp A_2x$ and $y \perp A_1x$. If such $y$ exists we have
	\begin{equation*}
y^*P(\lambda)x = (y^*A_2x)(\lambda^2 + \alpha_0)
	\end{equation*}
and condition \eqref{nozeroinD} is satisfied.
If such vector $y$ does not exist we have $\{A_1x\}^{\perp} \subseteq \{A_2x\}^{\perp}$ and consequently $A_2x = c_1 A_1x$ for some constant $c_1 \in \mathbb{C}$. Then:
	\begin{equation*}
P(\lambda)x = (c_1\lambda^2 + \lambda + c_1\alpha_0 + \beta_0)A_1x
	\end{equation*}
and since the polynomial $\lambda \mapsto P(\lambda)$ is regular and has no eigenvalues in $D$ we have that $A_1x \neq 0$ and  $\lambda \mapsto c_1\lambda^2 + \lambda + c_1\alpha_0 + \beta_0$ has no roots in $D$. Setting $y = A_1x$  we have
	\begin{equation*}
y^*P(\lambda)x = (A_1x)^*(c_1\lambda^2 + \lambda + c_1\alpha_0 + \beta_0)A_1x = \|A_1x\|^2(c_1\lambda^2 + \lambda + c_1\alpha_0 + \beta_0)
	\end{equation*}
and \eqref{nozeroinD} is again satisfied.

Similarly, if the polynomial $\lambda + \beta_0$ is stable with respect to $D$, then we seek for a vector $y$ such that $y \perp A_2x$ and $y \not\perp A_1x$. If such $y$ exists, then
	\begin{equation*}
y^*P(\lambda)x = (y^*A_2x)(\lambda + \beta_0)
	\end{equation*}
and \eqref{nozeroinD} is satisfied.
If such the vector $y$ does not exist we have $\{A_2x\}^{\perp} \subseteq \{A_1x\}^{\perp}$ and consequently $A_1x= c_2 A_2x$ for some constant $c_2 \in \mathbb{C}$ and 
	\begin{equation*}
P(\lambda)x =(\lambda^2 + c_2\lambda + \alpha_0 + c_2\beta_0)A_2x\text{.}
	\end{equation*}
As before, the scalar polynomial $\lambda \mapsto \lambda^2 + c_2\lambda + \alpha_0 + c_2\beta_0$ is stable and $A_2x \neq 0$. Taking $y = A_2x$ we have \eqref{nozeroinD} satisfied.

 Assume that \eqref{0notinD1} holds.
If there exists a vector $y \in \mathbb{C}^n\setminus\{0\}$ such that $y \perp A_2x, \; y \not\perp A_1x$ and $y \perp A_0x$, then 
	\begin{equation*}
y^* P(\lambda)x = (y^* A_1x)\lambda
	\end{equation*}
and \eqref{nozeroinD} is satisfied, due to the assumption that $0 \notin D$.

Hence, further of the proof we assume that such vector $y$ does not exist, i.e.
$
\{A_2x, A_0x\}^{\perp} \subseteq \{A_1x\}^{\perp}$,  equivalently $A_1x\in \Span\set{A_2x, A_0x}$.
Thus, we have 
\begin{equation}\label{A_1x}
A_1x = \alpha_0 A_2x + \beta_0 A_0x
\end{equation}
for some $\alpha_0, \beta_0 \in \Comp$. 
Next, we can write
	\begin{align*}
P(\lambda)x = \lambda^2 A_2x + \lambda(\alpha_0 A_2x + \beta_0 A_0x) + A_0x = \lambda(\lambda + \alpha_0)A_2x + (\beta_0\lambda + 1)A_0x\text{.}
	\end{align*}
	
We show now that at least one of the scalar polynomials $\lambda \mapsto \lambda(\lambda + \alpha_0)$ or $\lambda \mapsto \beta_0\lambda + 1$ is stable with respect to $D$.	
By the assumption that the polynomial $(z_1, z_2) \mapsto z_1z_2 A_2 + z_2  A_1 + A_0$ is stable with respect to $D^2$ we have that the equation $(z_1z_2 A_2 + z_2  A_1 + A_0)x = 0$ has no solutions $(z_1, z_2) \in D^2$. Substituting  \eqref{A_1x} we obtain that the equation
$$
z_2\left(z_1+\alpha_0\right) A_2 x+ \left( \beta_0 z_2 +1 \right)A_0x = 0
$$
has no solutions $(z_1, z_2) \in D^2$. This implies, that at least one of the scalar polynomials 
$\lambda \mapsto \lambda(\lambda + \alpha_0)$ or $\lambda \mapsto \beta_0\lambda + 1$ has no roots in $D$ and the claim follows.
%

If the polynomial $\lambda \mapsto \lambda(\lambda + \alpha_0)$ is stable with respect to $D$, then similarly as before we seek for a vector $y$ such that $y \not\perp A_2x$ and $y \perp A_0x$. If such $y$ exists we have
	\begin{equation*}
y^*P(\lambda)x = (y^*A_2x)\lambda(\lambda + \alpha_0)
	\end{equation*}
and condition \eqref{nozeroinD} is satisfied.
If such the vector $y$ does not exist we have $A_2x = c_1 A_0x$ for some constant $c_1 \in \mathbb{C}$. Then:
	\begin{equation*}
P(\lambda)x = \Big[c_1 \lambda^2 + (c_1\alpha_0 + \beta_0)\lambda + 1\Big]A_0x
	\end{equation*}
and since the polynomial $\lambda \mapsto P(\lambda)$ is regular and has no eigenvalues in $D$ we have that $A_0x \neq 0$ and $c_1\lambda^2 + (c_1\alpha_0 + \beta_0)\lambda + 1$ has no roots in $D$. Setting $y = A_0x$  we have
	\begin{equation*}
y^*P(\lambda)x = (A_0x)^*\Big[c_2 \lambda^2 + (c_2\alpha_0 + \beta_0)\lambda + 1\Big]A_0x = \|A_0x\|^2\Big[c_2 \lambda^2 + (c_2\alpha_0 + \beta_0)\lambda + 1\Big] 
	\end{equation*}
and \eqref{nozeroinD} is again satisfied.

Similarly, if the polynomial $\lambda \mapsto \beta_0\lambda + 1$ is stable with respect to $D$, then we seek for a vector $y$ such that $y \perp A_2x$ and $y \not\perp A_0x$. If such $y$ exists then
	\begin{equation*}
y^*P(\lambda)x = (y^*A_0x)(\beta_0\lambda + 1)
	\end{equation*}
and \eqref{nozeroinD} is satisfied.
If such vector $y$ does not exist we have $A_0x = c_2 A_2x$ for some constant $c_2 \in \mathbb{C}$ and 
	\begin{equation*}
P(\lambda)x = \Bigl[\lambda^2 + (\alpha_0 + c_2\beta_0 )\lambda + c_2\Bigr]A_2x\text{.}
	\end{equation*}
As before, the scalar polynomial $\lambda \mapsto \lambda^2 + (\alpha_0 + c_2\beta_0)\lambda + c_2$ is stable and $A_2x \neq 0$. Taking $y = A_2x$ we have \eqref{nozeroinD} satisfied.

Finally, assume that \eqref{0notinD2} holds. If there exists a vector $y \in \mathbb{C}^n\setminus\{0\}$ such that $y \not\perp A_2x, \; y \perp A_1x$ and $y \perp A_0x$, then 
	\begin{equation*}
y^* P(\lambda)x = (y^*A_2x)\lambda^2
	\end{equation*}
and \eqref{nozeroinD} is satisfied, due to the assumption that $0\notin D$.

Hence, further we assume that such vector $y$ does not exist, i.e.
$
\{A_1x, A_0x\}^{\perp} \subseteq \{A_2x\}^{\perp}$, equivalently $A_2x \in \Span\set{A_1x, A_0x}$.
Thus, we have
\begin{equation}\label{A_2x}
A_2x = \alpha_0 A_1x + \beta_0 A_0x
\end{equation}
for some $\alpha_0, \beta_0 \in \Comp$. 
Next, we can write
	\begin{align*}
P(\lambda)x = (\alpha_0\lambda^2 + \lambda)A_1x + (\beta_0\lambda^2 + 1)A_0x = \lambda(\alpha_0\lambda + 1)A_1x + (\beta_0\lambda^2 + 1)A_0x\text{.}
	\end{align*}
	
We show now that at least one of the scalar polynomials $\lambda \mapsto \lambda(\alpha_0\lambda + 1)$ or $\lambda \mapsto \beta_0\lambda^2 + 1$ is stable with respect to $D$.	
By the assumption that the polynomial $(z_1, z_2) \mapsto z_1^2z_2 A_2 + z_1^2 A_1 + z_2 A_0$ is stable with respect to $D^2$, we have that the equation $(z_1^2z_2 A_2 + z_1^2 A_1 + z_2 A_0)x = 0$ has no solutions $(z_1, z_2) \in D^2$. Substituting \eqref{A_2x} we obtain that the equation
$$
z_1^2z_2(\alpha_0 A_1x + \beta_0 A_0x) + z_1^2 A_1x + z_2 A_0x = 0
$$
this means the equation
$$
\label{noso}z_1^2(\alpha_0z_2 + 1)A_1x + z_2(\beta_0z_1^2 + 1)A_0x = 0
$$
has no solutions $(z_1, z_2) \in D^2$. This implies, that at least one of the scalar polynomials 
$\lambda \mapsto \lambda(\alpha_0\lambda + 1)$ or $\lambda \mapsto \beta_0\lambda^2 + 1$ has no roots in $D$. Otherwise, there would exist $\lambda_1, \lambda_2 \in D$ such that $\lambda_1(\alpha_0\lambda_1 + 1) = \beta_0\lambda_2^2 + 1 = 0$. Since $\lambda_1 \neq 0$ we would have $\alpha_0\lambda_1 + 1 = \beta_0\lambda_2^2 + 1 = 0$ and taking $(z_1, z_2) = (\lambda_2, \lambda_1) \in D^2$, which is a solution of the equation \eqref{noso} - contradiction. Therefore, the claim follows.
%

If the polynomial $\lambda \mapsto \lambda(\alpha_0\lambda + 1)$ is stable with respect to $D$, then similarly as before we seek for a vector $y$ such that $y \not\perp A_1x$ and $y \perp A_0x$. If such $y$ exists we have
	\begin{equation*}
y^*P(\lambda)x = (y^*A_1x)\lambda(\alpha_0\lambda + 1)
	\end{equation*}
and condition \eqref{nozeroinD} is satisfied.
If such vector $y$ does not exist we have $\{A_0x\}^{\perp} \subseteq \{A_1x\}^{\perp}$ and consequently $A_1x = c_1 A_0x$ for some constant $c_1 \in \mathbb{C}$. Then:
	\begin{equation*}
P(\lambda)x = \Bigl[\lambda(\alpha_0\lambda + 1)c_1 + \beta_0\lambda^2 + 1\Bigr]A_0x = \Bigl[(\alpha_0c_1 + \beta_0)\lambda^2 + c_1\lambda + 1\Bigr]A_0x
	\end{equation*}
and since the polynomial $\lambda \mapsto P(\lambda)$ is regular and has no eigenvalues in $D$ we have that $A_0x \neq 0$ and the polynomial $\lambda \mapsto (\alpha_0c_1 + \beta_0)\lambda^2 + c_1\lambda + 1$ has no roots in $D$. Setting $y = A_0x$  we have
	\begin{equation*}
y^*P(\lambda)x = (A_0x)^*\Bigl[(\alpha_0c_1 + \beta_0)\lambda^2 + c_1\lambda + 1\Bigr]A_0x = \|A_0x\|^2\Bigl[(\alpha_0c_1 + \beta_0)\lambda^2 + c_1\lambda + 1\Bigr]
	\end{equation*}
and \eqref{nozeroinD} is satisfied again.

Similarly, if the polynomial $\lambda \mapsto \beta_0\lambda^2 + 1$ is stable with respect to $D$, then we seek for a vector $y$ such that $y \perp A_1x$ and $y \not\perp A_0x$. If such $y$ exists, then
	\begin{equation*}
y^*P(\lambda)x = (y^*A_0x)(\beta_0\lambda^2 + 1)
	\end{equation*}
and \eqref{nozeroinD} is satisfied.
If such the vector $y$ does not exist we have $\{A_1x\}^{\perp} \subseteq \{A_0x\}^{\perp}$ and consequently $A_0x= c_2 A_1x$ for some constant $c_2 \in \mathbb{C}$ and 
	\begin{equation*}
P(\lambda)x = \Bigl[\lambda(\alpha_0\lambda + 1) + (\beta_0\lambda^2 + 1)c_2\Bigr]A_1x = \Bigl[(\alpha_0 + \beta_0c_2)\lambda^2 + \lambda + c_2\Bigr]A_1x\text{.}
	\end{equation*}
As before, the scalar polynomial $\lambda \mapsto (\alpha_0 + \beta_0c_2)\lambda^2 + \lambda + c_2$ is stable with respect to $D$ and $A_1x \neq 0$. Taking $y = A_1x$ we have \eqref{nozeroinD} satisfied.
	\end{proof}

From Theorem 3.6 (ii) we know that hyperstability of palindromic matrix polynomials of degree three is equivalent to their stability. In addition to this, we consider now polynomials of the form $P(\lambda) = \lambda^3 A_0 + \lambda^2 A_2 + \lambda A_1 + A_0$. Therefore, the following Theorem delivers us a lot of interesting examples, when $A_1 \neq A_2$.
    \begin{theorem}\label{poly3}
Let $P(\lambda) = \lambda^3 A_0 + \lambda^2 A_2 + \lambda A_1 + A_0$ be a cubic matrix polynomial and let $D$ be a nonempty open or closed subset of the complex plane $\mathbb{C}$. If at least one of the following conditions hold:
        \begin{enumerate}[\rm (a)]
\item\label{forA_0}
the multivariate matrix polynomial $(z_1, z_2) \mapsto (z_1^3 z_2^3 + z_1^3 + z_2^3)A_0 + (z_1^2 z_2^3 +z_1^2)A_2 + (z_1^3z_2 + z_2) A_1 + A_0$ is stable with respect to $D^2$ and $-1, \frac{1}{2} - \frac{\sqrt{3}}{2}i, \frac{1}{2} + \frac{\sqrt{3}}{2}i \not\in D$,
\item\label{forA_1}
the multivariate matrix polynomial $(z_1, z_2) \mapsto z_2^3 A_0 + z_1 z_2 A_2 + z_2 A_1 + A_0$ is stable with respect to $D^2$ and $0 \not\in D$,
\item\label{forA_2}
the multivariate matrix polynomial $(z_1, z_2) \mapsto z_1 z_2^3 A_0 + z_1 z_2^2 A_2 + z_2^2 A_1 + z_1 A_0$ is stable with respect to $D^2$ and $0 \not\in D$,
        \end{enumerate}
then the matrix polynomial $ P(\lambda)$ is hyperstable with respect to $D$.
    \end{theorem}
    \begin{proof}
In this proof, we proceed as in the proof of Theorem 5.2. Setting $z_1 = z_2 = \lambda$ leads in each case to stability of the polynomial $\lambda \mapsto P(\lambda)$ with respect to $D$. In particular, $P(\lambda) $ is regular.  In the case \eqref{forA_0} stability of a polynomial $\lambda \mapsto (\lambda^3 + 1)P(\lambda)$ is equivalent to stability of the polynomial $\lambda \mapsto P(\lambda)$. For hyperstability of the polynomial $\lambda \mapsto P(\lambda)$ with respect to $D$, fix $x \in \mathbb{C}^n\setminus\{0\}$.

Assume that the condition \eqref{forA_0} holds. If there exists a vector $y \in \mathbb{C}^n\setminus\{0\}$ such that $y \not\perp A_0x, y \perp A_2x, y \perp A_1x$, then
    \begin{equation}
y^*P(\lambda)x = (y^*A_0x)(\lambda^3 + 1)
    \end{equation}
and the condition \eqref{nozeroinD} is satisfied because of fact that $-1, \frac{1}{2} - \frac{\sqrt{3}}{2}i, \frac{1}{2} + \frac{\sqrt{3}}{2}i \not\in D$. Otherwise, we have $\{A_2x, A_1x\}^{\perp} \subseteq \{A_0x\}^{\perp}$, which equivalently means that $A_0x = \alpha_0 A_2x + \beta_0 A_1x$ for some $\alpha_0, \beta_0 \in \mathbb{C}$. Thus we can write
    \begin{equation*}
P(\lambda)x = (\alpha_0\lambda^3 + \lambda^2 + \alpha_0)A_2x + (\beta_0\lambda^3 + \lambda + \beta_0)A_1x\text{.}
    \end{equation*}
Due to stability of the multivariate polynomial in \eqref{forA_0} we conclude that the equation
    \begin{equation*}
\Bigl[(z_1^3 z_2^3 + z_1^3 + z_2^3)A_0 + (z_1^2 z_2^3 +z_1^2)A_2 + (z_1^3z_2 + z_2) A_1 + A_0\Bigr]x = 0
    \end{equation*}
does not have solutions $(z_1, z_2) \in D^2$. Then the equation
    \begin{equation*}
(\alpha_0 z_1^3 z_2^3 + z_1^2 z_2^3 + \alpha_0 z_1^3 + \alpha_0 z_2^3 + z_1^2 + \alpha_0)A_2x + (\beta_0 z_1^3 z_2^3 + z_1^3 z_2 + \beta_0 z_1^3 + \beta_0 z_2^3 + z_2 + \beta_0)A_1x = 0
    \end{equation*}
does not have such solutions as well. Therefore, at least one of the polynomials $\lambda \mapsto \alpha_0\lambda^3 + \lambda^2 + \alpha_0$ or $\lambda \mapsto \beta_0\lambda^3 + \lambda + \beta_0$ is stable with respect to $D$. If the polynomial $\lambda \mapsto \alpha_0\lambda^3 + \lambda^2 + \alpha_0$ is stable with respect $D$, then we seek for a vector $y \not\perp A_2x$ and $y \perp A_1x$. If such vector $y$ exists, then
    \begin{equation*}
y^*P(\lambda)x = (y^*A_2x)(\alpha_0\lambda^3 + \lambda^2 + \alpha_0)
    \end{equation*}
and the condition \eqref{nozeroinD} is also satisfied. Otherwise, we have $\{A_1x\}^{\perp} \subseteq \{A_2x\}^{\perp}$, which equivalently means that $A_2x = c_1A_1x$ for some $c_1 \in \mathbb{C}$. Finally, we can write
    \begin{equation*}
P(\lambda)x = \Bigl[(c_1\alpha_0 + \beta_0)\lambda^3 + x_1\lambda^2 + \lambda + (c_1\alpha_0 + \beta_0)\Bigr]A_1x,
    \end{equation*}
where a polynomial in the square brackets is stable because of stability of the polynomial $\lambda \mapsto P(\lambda)$. Now, we just take $y = A_1x$ to satisfy the condition \eqref{nozeroinD}:
    \begin{equation*}
y^*P(\mu)x = ||A_1x||^2\Bigl[(c_1\alpha_0 + \beta_0)\mu^3 + c_1\mu^2 + \mu + (c_1\alpha_0 + \beta_0)\Bigr] \neq 0
    \end{equation*}
for $\mu \in D$. If the polynomial $\lambda \mapsto \beta_0\lambda^3 + \lambda + \beta_0$ is stable with respect $D$, then we seek in turn for a vector $y \perp A_2x$ and $y \not\perp A_1x$. If such vector $y$ exists, then
    \begin{equation*}
y^*P(\lambda)x = (y^*A_1x)(\beta_0\lambda^3 + \lambda + \beta_0)
    \end{equation*}
and the condition \eqref{nozeroinD} is satisfied again. Otherwise, we have $\{A_2x\}^{\perp} \subseteq \{A_1x\}^{\perp}$, which equivalently means that $A_1x = c_2A_2x$ for some $c_2 \in \mathbb{C}$. As before, we can write:
\begin{equation*}
P(\lambda)x = \Bigl[(\alpha_0 + c_2\beta_0)\lambda^3 + \lambda^2 + c_2\lambda + (\alpha_0 + c_2\beta_0)\Bigr]A_2x,
    \end{equation*}
where a polynomial in the square brackets is stable because of stability of the polynomial $\lambda \mapsto P(\lambda)$. As previously, we take $y = A_2x$ to satisfy the condition \eqref{nozeroinD}:
    \begin{equation*}
y^*P(\mu)x = ||A_2x||^2\Bigl[(\alpha_0 + c_2\beta_0)\lambda^3 + \lambda^2 + c_2\lambda + (\alpha_0 + c_2\beta_0)\Bigr] \neq 0
    \end{equation*}
for $\mu \in D$, which ends the proof in this case.

Assume that the condition \eqref{forA_1} holds. If there exists a vector $y \in \mathbb{C}^n\setminus\{0\}$ such that $y \perp A_0x, y \perp A_2x, y \not\perp A_1x$, then
    \begin{equation}
y^*P(\lambda)x = (y^*A_1x)\lambda
    \end{equation}
and the condition\eqref{nozeroinD} is satisfied because of fact that $0 \not\in D$. Otherwise, we have $\{A_0x, A_2x\}^{\perp} \subseteq \{A_1x\}^{\perp}$, which equivalently means that $A_1x = \alpha_0 A_2x + \beta_0 A_0x$ for some $\alpha_0, \beta_0 \in \mathbb{C}$. Thus we can write
    \begin{equation*}
P(\lambda)x = (\lambda^2 + \alpha_0\lambda)A_2x + (\lambda^3 +\beta_0\lambda +1)A_0x\text{.}
    \end{equation*}
Due to stability of the multivariate polynomial in \eqref{forA_1} we conclude that the equation
    \begin{equation*}
\Bigl[z_2^3 A_0 + z_1 z_2 A_2 + z_2 A_1 + A_0\Bigr]x = 0
    \end{equation*}
does not have solutions $(z_1, z_2) \in D^2$. Then the equation
    \begin{equation*}
(z_1 z_2 + \alpha_0 z_2)A_2x + (z_2^3 + \beta_0z_2 + 1)A_0x = 0
    \end{equation*}
does not have such solutions as well. Therefore, at least one of the polynomials $\lambda \mapsto \lambda^2 + \alpha_0\lambda$ or $\lambda \mapsto \lambda^3 + \beta_0\lambda + 1$ is stable with respect to $D$. If the polynomial $\lambda \mapsto \lambda^2 + \alpha_0\lambda$ is stable with respect $D$, then we seek for a vector $y \not\perp A_2x$ and $y \perp A_0x$. If such vector $y$ exists, then
    \begin{equation*}
y^*P(\lambda)x = (y^*A_2x)(\lambda^2 + \alpha_0\lambda)
    \end{equation*}
and the condition \eqref{nozeroinD} is also satisfied. Otherwise, we have $\{A_0x\}^{\perp} \subseteq \{A_2x\}^{\perp}$, which equivalently means that $A_2x = c_1A_0x$ for some $c_1 \in \mathbb{C}$. Finally, we can write
    \begin{equation*}
P(\lambda)x = \Bigl[\lambda^3 + c_1\lambda^2 + (c_1\alpha_0 + \beta_0)\lambda + 1\Bigr]A_0x,
    \end{equation*}
where a polynomial in the square brackets is stable because of stability of the polynomial $\lambda \mapsto P(\lambda)$. Now, we just take $y = A_0x$ to satisfy the condition \eqref{nozeroinD}:
    \begin{equation*}
y^*P(\mu)x = ||A_0x||^2\Bigl[\lambda^3 + c_1\lambda^2 + (c_1\alpha_0 + \beta_0)\lambda + 1\Bigr] \neq 0
    \end{equation*}
for $\mu \in D$. If the polynomial $\lambda \mapsto \lambda^3 + \beta_0\lambda + 1$ is stable with respect $D$, then we seek in turn for a vector $y \perp A_2x$ and $y \not\perp A_0x$. If such vector $y$ exists, then
    \begin{equation*}
y^*P(\lambda)x = (y^*A_0x)(\lambda^3 + \beta_0\lambda + 1)
    \end{equation*}
and the condition \eqref{nozeroinD} is satisfied again. Otherwise, we have $\{A_2x\}^{\perp} \subseteq \{A_0x\}^{\perp}$, which equivalently means that $A_0x = c_2A_2x$ for some $c_2 \in \mathbb{C}$. As before, we can write:
\begin{equation*}
P(\lambda)x = \Bigl[c_2\lambda^3 + \lambda^2 + (\alpha_0 + c_2\beta_0)\lambda + c_2\Bigr]A_2x,
    \end{equation*}
where a polynomial in the square brackets is stable because of stability of the polynomial $\lambda \mapsto P(\lambda)$. As previously, we take $y = A_2x$ to satisfy the condition \eqref{nozeroinD}:
    \begin{equation*}
y^*P(\mu)x = ||A_2x||^2\Bigl[c_2\lambda^3 + \lambda^2 + (\alpha_0 + c_2\beta_0)\lambda + c_2\Bigr] \neq 0
    \end{equation*}
for $\mu \in D$, which ends the proof in this case.

Assume that the condition \eqref{forA_2} holds. If there exists a vector $y \in \mathbb{C}^n\setminus\{0\}$ such that $y \perp A_0x, y \not\perp A_2x, y \perp A_1x$, then
    \begin{equation}
y^*P(\lambda)x = (y^*A_2x)\lambda^2
    \end{equation}
and the condition\eqref{nozeroinD} is satisfied because of fact that $0 \not\in D$. Otherwise, we have $\{A_1x, A_0x\}^{\perp} \subseteq \{A_2x\}^{\perp}$, which equivalently means that $A_2x = \alpha_0 A_1x + \beta_0 A_0x$ for some $\alpha_0, \beta_0 \in \mathbb{C}$. Thus we can write
    \begin{equation*}
P(\lambda)x = (\alpha_0\lambda^2 + \lambda)A_2x + (\lambda^3 + \beta_0\lambda^2 + 1)A_0x\text{.}
    \end{equation*}
Due to stability of the multivariate polynomial in \eqref{forA_2} we conclude that the equation
    \begin{equation*}
\Bigl[z_1 z_2^3 A_0 + z_1 z_2^2 A_2 + z_2^2 A_1 + z_1 A_0\Bigr]x = 0
    \end{equation*}
does not have solutions $(z_1, z_2) \in D^2$. Then the equation
    \begin{equation*}
(\alpha_0 z_1 z_2^2 + z_2^2)A_1x + (z_1 z_2^3 + \beta_0 z_1 z_2^2 + z_1)A_0x = 0
    \end{equation*}
does not have such solutions as well. Therefore, at least one of the polynomials $\lambda \mapsto \alpha_0\lambda^2 + \lambda$ or $\lambda \mapsto \lambda^3 + \beta_0\lambda^2 + 1$ is stable with respect to $D$. If the polynomial $\lambda \mapsto \alpha_0\lambda^2 + \lambda$ is stable with respect $D$, then we seek for a vector $y \not\perp A_1x$ and $y \perp A_0x$. If such vector $y$ exists, then
    \begin{equation*}
y^*P(\lambda)x = (y^*A_1)(\alpha_0\lambda^2 + \lambda)
    \end{equation*}
and the condition \eqref{nozeroinD} is also satisfied. Otherwise, we have $\{A_0x\}^{\perp} \subseteq \{A_1x\}^{\perp}$, which equivalently means that $A_1x = c_1A_0x$ for some $c_1 \in \mathbb{C}$. Finally, we can write
    \begin{equation*}
P(\lambda)x = \Bigl[\lambda^3 + (c_1\alpha_0 + \beta_0)\lambda^2 + c_1\lambda + 1\Bigr]A_0x,
    \end{equation*}
where a polynomial in the square brackets is stable because of stability of the polynomial $\lambda \mapsto P(\lambda)$. Now, we just take $y = A_0x$ to satisfy the condition \eqref{nozeroinD}:
    \begin{equation*}
y^*P(\mu)x = ||A_0x||^2\Bigl[\lambda^3 + (c_1\alpha_0 + \beta_0)\lambda^2 + c_1\lambda + 1\Bigr] \neq 0
    \end{equation*}
for $\mu \in D$. If the polynomial $\lambda \mapsto \lambda^3 + \beta_0\lambda^2 + 1$ is stable with respect $D$, then we seek in turn for a vector $y \perp A_1x$ and $y \not\perp A_0x$. If such vector $y$ exists, then
    \begin{equation*}
y^*P(\lambda)x = (y^*A_0x)(\lambda^3 + \beta_0\lambda^2 + 1)
    \end{equation*}
and the condition \eqref{nozeroinD} is satisfied again. Otherwise, we have $\{A_1x\}^{\perp} \subseteq \{A_0x\}^{\perp}$, which equivalently means that $A_0x = c_2A_1x$ for some $c_2 \in \mathbb{C}$. As before, we can write:
\begin{equation*}
P(\lambda)x = \Bigl[c_2\lambda^3 + (\alpha_0 + c_2\beta_0)\lambda^2 + \lambda + c_2\Bigr]A_1x,
    \end{equation*}
where a polynomial in the square brackets is stable because of stability of the polynomial $\lambda \mapsto P(\lambda)$. As previously, we take $y = A_1x$ to satisfy the condition \eqref{nozeroinD}:
    \begin{equation*}
y^*P(\mu)x = ||A_1x||^2\Bigl[c_2\lambda^3 + (\alpha_0 + c_2\beta_0)\lambda^2 + \lambda + c_2\Bigr] \neq 0
    \end{equation*}
for $\mu \in D$, which ends the proof in this case and the whole proof is completed.
\end{proof}

\section{The polarization operator for  matrix polynomials}\label{sPol}

In this Section we develop the theory  of hyperstable matrix polynomials of several variables. 
This, besides an independent interest,  will later on  serve, once again, as a tool for showing hyperstability of matrix polynomial of one variable.

 First, let us define a scalar multi-variate polynomial $(z_1, z_2, \dots , z_{\kappa}) \mapsto p(z_1, z_2, \dots, z_\kappa)$  to be \emph{multi-affine} if and only if its degree with respect to each variable $z_j$, $j=1,\dots,\kappa$,  is less or equal to one. Further, we say that $(z_1, z_2, \dots , z_{\kappa}) \mapsto p(z_1, z_2, \dots, z_\kappa)$ is \emph{symmetric} if and only if any permutation of the variables $z_1, z_2, \dots, z_\kappa$ leaves the polynomial intact.  
Below we recall the famous Grace-Walsh-Szeg\"o Theorem \cite{Gr1902,Sze1922,Wa1922}, we present its shortened  version needed for the current investigations.

\begin{theorem}\label{GWS}
Let $p \in \mathbb{C}[z_1, z_2, \dots , z_{\kappa}]$ be a symmetric multi-affine polynomial, let $D$ be an open or closed disc or open or closed half-plane  and let $\zeta_1, \zeta_2, \dots , \zeta_{\kappa}\in D$. Then there exists a point $\zeta_0 \in D$ such that
		\begin{equation*}
p(\zeta_1, \zeta_2, \dots , \zeta_{\kappa}) = p(\zeta_0, \zeta_0, \dots , \zeta_0)\text{.}
		\end{equation*}
	\end{theorem}
	We introduce now  the main object of this Section.  Let $\kappa \in \mathbb{Z}_+$ and let $P(\lambda)=\lambda^d A_d + \lambda^{d-1} A_{d-1} + \dots + A_0$ be a matrix polynomial of degree $d \leq \kappa$, we put $A_{d+1} = A_{d+2} = \dots = A_{\kappa} = 0$. We define the \emph{polarization} of $P(\lambda)$ by 
		\begin{equation}\label{Tdef1}
(T_{\kappa}P) (z_1, z_2, \dots , z_{\kappa}):= \sum_{j=0}^{\kappa} \binom{\kappa}{j}^{-1}s_j(z_1, z_2, \dots , z_{\kappa})A_j, 
		\end{equation}
where the symbol $s_j$ denotes the $j$-th elementary symmetric polynomial
		\begin{equation}\label{Tdef2}
s_0(z_1, z_2, \dots, z_{\kappa}) := 1,\;\;\; s_j(z_1, z_2, \dots , z_{\kappa}) := \sum_{1 \leq i_1 < i_2 < \dots < i_j \leq \kappa} z_{i_1}z_{i_2} \dots z_{i_j}.
		\end{equation}
The operator $T_\kappa$ defined above is a well known tool for scalar multi-variate polynomials, see, e.g., \cite{BorB09}. Our contribution is to extend its action to matrix polynomials. Note that the operator $T_\kappa$ is injective with its image being the set consisting of all symmetric multi-affine polynomials (with matrix coefficients). We have the following result. 	



	\begin{theorem}\label{Tkappa2}
 Let  $\kappa\geq 1$ and let the operator $T_\kappa$ be defined by \eqref{Tdef1} above with $D$ being an open or closed disc or open or closed half-plane. Then a matrix polynomial $P(\lambda)$ is hyperstable with respect to  $D$ if and only if $(T_\kappa P)(z_1,\dots,z_\kappa)$  is hyperstable with respect to $D^\kappa$.
	\end{theorem}	

	\begin{proof}
Consider first the case $n=1$, i.e., take a scalar polynomial $p(\lambda)$ and suppose that there exists a point $(\zeta_1, \zeta_2, \dots , \zeta_{\kappa}) \in D^{\kappa}$ such that $(T_{\kappa}p)(\zeta_1, \zeta_2, \dots , \zeta_{\kappa}) = 0$. Since the polynomial $(z_1, z_2, \dots , z_{\kappa}) \mapsto (T_\kappa p)(z_1, z_2, \dots z_{\kappa})$ is a symmetric multi-affine polynomial, then from Theorem~\ref{GWS} we conclude that there exists  a point $\zeta_0 \in D$ such that $(T_\kappa p)(\zeta_1, \zeta_2, \dots , \zeta_{\kappa}) = (T_\kappa p)(\zeta_0, \zeta_0, \dots , \zeta_0)$. However, note that
$ (T_\kappa p)(\zeta_0, \zeta_0, \dots , \zeta_0)=p(\zeta_0)$, which shows the forward implication.
The converse implication is obvious, as $(T_{\kappa}p)(\lambda, \lambda, \dots , \lambda) = p(\lambda)$.	
	
Assume now that $P(\lambda)$ is a matrix polynomial, hyperstable with respect to $D$. Take an arbitrary  vector $x \in \mathbb{C}^n \setminus\set0$. By definition of hyperstability, there exists $y\in \mathbb{C}^n \setminus\set0$ such that the scalar polynomial $p(\lambda)=y^*P(\lambda)x$ is stable with respect to $D$. By the first part of the proof,  the polynomial $(T_\kappa p)(z_1,\dots,z_\kappa)$ is stable with respect to $D^\kappa$. Observe that 
$$
(T_\kappa p)(z_1,\dots,z_\kappa)= y^* \left( (T_\kappa P   )(z_1,\dots,z_\kappa)\right)x,
$$
 which shows that  $(T_\kappa P   )(z_1,\dots,z_\kappa)$ is hyperstable with respect to $D^\kappa$.  
 The converse implication is again obvious.
\end{proof}

The first part of the proof is well known (in a slightly different setting), see e.g. \cite{BorB09}, and is presented here for completeness.
One may wonder if the following version of Theorem~\ref{Tkappa2} is true: if $P(\lambda)$ is stable with respect to $D$ then $(T_\kappa P)(z_1,\dots,z_\kappa)$ is stable with respect to $D$. However, it is not the case.

\begin{example}\label{nonstab}\rm  We continue with $P(\lambda)$ as in Example \ref{exa}. Again, $P(\lambda)$ is stable with respect to any $D$ but 
$$
\det((T_2 P)(z_1,z_2))=\det\mat{cc} 1 & \frac{z_1+z_2}2 \\ \frac{z_1+z_2}2 & z_1z_2+1\rix=1-\left( \frac{z_1-z_2}2\right)^2. 
$$
Hence, every $(\mu_1,\mu_2)$ with $\mu_1-\mu_2=\pm 2$ is an eigenvalue, in particular $(T_2P)(z_1,z_2)$ is not stable with respect to, e.g., $H_0^2$. 
\end{example}

We present now a  general tool for creating hyperstable  matrix polynomials of one variable using the operator $T_\kappa$. Its applications will be given in the next Section. 

\begin{theorem}\label{increasedegree}
Let $D$ be an open or closed disc or open or closed half-plane. If a matrix polynomial $P(\lambda)$ of degree $d $ is hyperstable with respect to $D$, then for any scalar polynomials $p_1(\lambda), p_2(\lambda) \dots, p_\kappa(\lambda)$, $\kappa\geq d$, the matrix polynomial $Q(\lambda) := (T_{\kappa}P)(p_1(\lambda), p_2(\lambda), \dots , p_\kappa(\lambda))$ is hyperstable with respect to
$E := p_1^{-1}(D) \cap p_2^{-1}(D) \cap \dots \cap p_\kappa^{-1}(D)\subseteq\Comp$.
\end{theorem}
	
\begin{proof}
Fix a nonzero vector $x\in\Comp^n\setminus\{0\}$. By Theorem~\ref{Tkappa2}, the multivariate polynomial $(T_{\kappa}P)(z_1, z_2, \dots , z_\kappa)$ is hyperstable with respect to $D^\kappa$, i.e., there exists $y\in\Comp^n\setminus\{0\}$ such that $y^*(T_{\kappa}P)(z_1, z_2, \dots , z_\kappa)x \neq 0$ for all $z_1, z_2, \dots, z_\kappa \in D$. 
In particular, for any $\lambda\in\Comp$ such that $p_j(\lambda)\in D$ ($j=1,\dots,\kappa$) one has 
	$y^*(T_{\kappa}P)(p_1(\lambda), p_2(\lambda), \dots , p_\kappa(\lambda))x \neq 0$. Hence, $Q(\lambda)$ is hyperstable with respect to $E$.
\end{proof}
	
\section{Hyperstability and stability of some classes of matrix polynomials}\label{sPos}

Now we will apply the theory developed in previous Sections. First let us deal with a relatively simple matrix polynomial, connected with the  Moore-Gibson-Thompson eigenvalue equation \cite{Ben21,KalN19}.
	\begin{proposition}\label{MGT} 
	Let $P(\lambda)=\lambda^3 I_n +a I_n \lambda^2 +\lambda b R + cR$ with  $R\in\Comp^{n,n}$ positive definite. If $a>1$ and $b>c$ then $P(\lambda)$ is hyperstable with respect to the open right half-plane.
\end{proposition}

\begin{proof}
	It was showed in  \cite{MehMW22}  that $P(\lambda)$ is stable with respect to the open right half-plane. By Theorem~\ref{uppert}\eqref{pq} it is also  hyperstable. 
\end{proof}

Second, let us explore a simple triangle inequality argument, which shows an application of Theorem~\ref{poly2}. Note that stability of $P(\lambda)$ below  is almost obvious, but  hyperstability needs some justification.

\begin{proposition}\label{subadd}
 Let $A_0, A_1, A_2 \in \mathbb{C}^{n \times n}$ and let $D = \set{z \in \Comp: |z|<r}$, $r>0$ be such that  $r\norm{A_1}+r^2\norm{A_2}<\sigma_{\min}(A_0)$. Then the multivariate polynomial $z_1^2A_2+z_2A_1+A_0$ is stable with respect to $D$ and consequently 
 $P(\lambda) = \lambda^2 A_2 + \lambda A_1 + A_0$ is a regular matrix polynomial, hyperstable with respect to $D$.
 \end{proposition}

\begin{proof}
 Note that
$$
\norm{z_1^2 A_2+z_2A_1}\leq r\norm {A_1}+r^2\norm{A_2}<\sigma_{\min}(A_0),
$$
which implies stability of the multivariate polynomial $z_1^2A_2+z_2A_1+A_0$ with respect to $D$.
Application of Theorem~\ref{poly2}\eqref{0?D} finishes the proof. 
\end{proof}


The quadratic polynomials with $A_2$, $A_0$  Hermitian positive semi-definite and the Hermitian part of $ A_1$  positive semi-definite were studied extensively in \cite{MehMW18,MehMW21}. It was shown therein that they are stable with respect to the right half-plane. We show more, namely hyperstablity with respect to the right half-plane.

\begin{theorem}\label{half-plane} Let $R_j \in \mathbb{C}^{n, n}$ $(j=0,1,2)$  be  Hermitian positive semi-definite and let $J\in \mathbb{C}^{n, n}$ be skew-Hermitian. Assume that $\ker R_0\cap\ker R_1\cap \ker R_2\cap\ker J=\set0$.
 Then the polynomial $z_1z_2 R_2 +z_2 (J+R_1) +R_0$ is stable with respect to  $H_{\pi/2}^2$. In consequence,  $P(\lambda) = \lambda^2 R_2 + \lambda (J+R_1) + R_0$ is a regular polynomial, hyperstable  with respect to the open right half-plane $H_{\pi/2}$.
 \end{theorem}
\begin{proof}
We will make use of Theorem \ref{poly2} (b), note that $0 \not\in D = H_0$. 
Consider the polynomial $\tilde P(z_1,z_2)=z_1z_2R_2+z_2 (J+R_1) +R_0$, and suppose it is not stable, i.e., $\tilde P(\mu_1,\mu_2)x=0$ for some  $(\mu_1,\mu_2) \in H_{\pi/2}^2$ and $x\neq 0$. Multiplying from the left by $x^*$ and taking the real part one obtains
\begin{equation}\label{eigenvalue1}
\re(\mu_1) x^*R_2 x +  x^*R_1x + \re\left(\frac1{\mu_2}\right) x^*R_0 x=0.
\end{equation}
 Since both $\re( \mu_1)$ and $\re(\frac1{\mu_2})$ are positive and $R_2,R_1,R_0$ are positive semi-definite, we obtain  $x^*R_2 x =  x^*R_1x  = x^*R_0 x=0$. Hence, $R_2x=R_1x=R_0x=0$.
 But this implies that $0=\tilde P(\mu_1,\mu_2)x=Jx$, 
 contradiction.

	
	\end{proof}
	
\begin{remark}\rm
    Let us remark two things. First, the condition $\ker R_1\cap\ker R_2\cap \ker R_3\cap\ker J=\set0$ is equivalent to regularity of $P(\lambda)$, see \cite{MehMW21}. Second, although we multiplied in the proof from the left by $x^*$, we have \emph{not} shown that the numerical range of $P(\lambda)$ is outside the open right half-plane.   
\end{remark}	
	
\begin{corollary}\label{c-hp}
 Let $R\in\mathbb C^{n,n}$  be  Hermitian positive semi-definite,   let $J\in \mathbb{C}^{n, n}$ be skew-Hermitian and let  $Q, A_0, A_2 \in \mathbb{C}^{n, n}$ be such that $Q^*A_2$ and $Q^*A_0$ are Hermitian positive semi-definite.  Assume also  that $\ker Q^* A_0\cap\ker(Q^* R Q)\cap \ker (Q^*JQ)\cap\ker Q^*A_2=\set0$.
 Then the matrix polynomial  $P(\lambda) = \lambda^2 A_2 + \lambda (J+R)Q + A_0$ is a regular polynomial, hyperstable  with respect to the open right half-plane $H_{\pi/2}$.
\end{corollary}	
	
	\begin{proof}
	First, observe that the polynomial $ \lambda^2 Q^*A_2 + \lambda Q^*(J+R)Q + Q^*A_0$ satisfies the assumptions of Theorem~\ref{half-plane}. Second, apply Lemma~\ref{lQ}.
	\end{proof}
	
In what follows,  the argument of a complex number is taken such that $\Arg z \in (-\pi; \pi]$ and for the sake of simplicity we set $\Arg 0 := 0$.
	It was shown in~\cite{MehMW22} that a  regular cubic matrix polynomial with all coefficients positive semi-definite is stable with respect to $D = \{z \in \mathbb{C} : -\pi/3 < \Arg z < \pi/3\}$. We present here a related result.
	    \begin{theorem}\label{ker}
Let $R_j \in \mathbb{C}^{n, n}$ $(j=1,2,3)$  be  Hermitian positive semi-definite and let $A_0\in\mathbb C^{n,n}$ be Hermitian and let $G\in \mathbb{C}^{n, n}$ be skew-Hermitian with $-\ii G$ positive semi-definite. Assume that $\ker G \cap \ker A_0\cap \ker R_1\cap\ker R_2\cap \ker R_3 = \set0$.
 Then the following holds:
        \begin{enumerate}[\rm (i)]
\item  the multivariate matrix polynomial $P_1(z_1, z_2)= (z_1^3 z_2^3 + z_1^3 + z_2^3)R_3 + (z_1^2 z_2^3 +z_1^2)R_2 + (z_1^3z_2 + z_2) R_1 + A_0+G$ is stable with respect to $D_1^2$, where $D_1 = \{z \in \mathbb{C} : 0 < \Arg z < \pi/6\}$;
\item\label{P2} the multivariate matrix polynomial $P_2(z_1, z_2) = z_2^3 R_3 + z_1 z_2 R_2 + z_2 R_1 + A_0+G$ is stable with respect to $D_2^2$, where $D_2 = \{z \in \mathbb{C} : 0 < \Arg z < \pi/3\}$;
\item the multivariate matrix polynomial $
P_3(z_1, z_2) = z_1 z_2^3 R_3 + z_1 z_2^2 R_2 + z_2^2 R_1 +  z_1 (A_0 + G)$ is stable with respect to $D_3^2$, where  $D_3 = \{z \in \mathbb{C} : 0 < \Arg z < \pi/4\}$.
        \end{enumerate}
In particular, $P(\lambda) = \lambda^3 R_3+\lambda^2 R_2+\lambda R_1 +A_0+G$ is stable with respect to $D = \{z \in \mathbb{C} : 0 < \Arg z < \pi/3\}$.
    \end{theorem}
    \begin{proof}
We show only (i), the proofs of (ii) and (iii) are very similar. Suppose that the matrix polynomial $(z_1, z_2) \mapsto P_1(z_1, z_2)$ has an eigenvalue $(\mu_1, \mu_2) \in D_1^2$. Thus there exist a nonzero vector $x \in \mathbb{C}^n\setminus\{0\}$ such that
        \begin{equation}\label{mmm}
(\mu_1^3\mu_2^3 + \mu_1^3 + \mu_2^3)R_3x + (\mu_1^2\mu_2^3 + \mu_1^2)R_2x + (\mu_1^3\mu_2 + \mu_2)R_1x + (A_0+G)x = 0\text{.}
        \end{equation}
Multiplying by $x^*$ and taking the imaginary part of both sides of the equation above we obtain
        \begin{equation*}
(x^*R_3x)\;\im(\mu_1^3\mu_2^3 + \mu_1^3 + \mu_2^3) + (x^*R_2x)\;\im(\mu_1^2\mu_2^3 + \mu_1^2) + (x^*R_1x)\;\im(\mu_1^3\mu_2 + \mu_2) -  x^*(\ii G)x = 0\text{.}
        \end{equation*}
Note that by assumption $x^*R_3x, x^*R_1x, x^*R_2x, -x^*(\ii G)x \geq 0$ and since $0< \Arg \mu_1, \Arg \mu_2 < \pi/6$, we have $\im(\mu_1^3\mu_2^3 + \mu_1^3 + \mu_2^3), \im(\mu_1^2\mu_2^3 + \mu_1^2), \im(\mu_1^3\mu_2 + \mu_2) > 0$. Therefore, we have $x^*R_3x = x^*R_2x = x^*R_1x = x^*(iG)x = 0$ and consequently  $R_3x = R_2x = R_1x = G x = 0$. Due to the equation \eqref{mmm} we have $A_0x = 0$, contradiction. 
The `In particular' part follows by substituting $\lambda$ for $z_1$ and $z_2$ in \eqref{P2}.
    \end{proof}
    Directly from Theorems~\ref{poly3} and \ref{ker}\eqref{P2} we get the following Corollary. 
    
    \begin{corollary} Let $R_j \in \mathbb{C}^{n, n}$ $(j=0,1,2)$  be  Hermitian positive semi-definite, then the polynomial 
 $P(\lambda) = \lambda^3 R_0 + \lambda^2 R_2 + \lambda  R_1 + R_0$ is hyperstable with respect to $ D = \{z \in \mathbb{C} : 0 < \Arg z < \pi/3\}$.
    \end{corollary}
We show now how the polarization operator may be used to increase the degree of the polynomial. The price for it is narrowing the set $D$.
	
\begin{corollary}\label{deg3} Let $R_j \in \mathbb{C}^{n, n}$ $(j=0,1,2)$  be  Hermitian positive semi-definite, and let $J\in\Comp^{n,n}$ be skew-Hermitian. Then the matrix polynomial 
  $Q(\lambda) = \lambda^3 R_2 + (\lambda^2+\lambda)(R_1+J)  + R_0$ is a regular cubic matrix polynomial,  hyperstable with respect to the angle $E := \{z\in\Comp:-\pi/4 < \Arg z  <\pi/4\}\setminus\{0\}$.
 \end{corollary}
	\begin{proof} By Theorem~\ref{half-plane} we get the polynomial $P(\lambda)= \lambda^2 R_2 + 2 \lambda  (R_1+J) + R_0$ hyperstable with respect to the open right half-plane $H_{\pi/2}$. We apply now Theorem~\ref{increasedegree} with $p_1(\lambda)=\lambda^2$, $p_2(\lambda)=\lambda$ getting 
	$$
	(T_2)P(z_1,z_2)=  z_1z_2 A_2 + (z_1+z_2) A_1 + A_0
	$$
	so that
		$
		(T_2P)(\lambda^2,\lambda) = Q(\lambda)$.
		Finally, observe that the angle  $E$ is precisely the set $p_1^{-1}(H_{\pi/2})\cap p_2^{-1}(H_{\pi/2})$ from Theorem~\ref{increasedegree}.
	\end{proof}
	
In Theorem~\ref{deg3} the operator $T_2$ was used to increase the degree of the quadratic polynomial $P(\lambda)$, which is known to by hyperstable. Let us show now an opposite action, where we use $T_2$ to decrease the degree.

\begin{corollary}
Let $R_j \in \mathbb{C}^{n, n}$ $(j=0,1)$  be  Hermitian positive semi-definite, and let $J\in\Comp^{n,n}$ be skew-Hermitian. Consider a  linear pencil $P(\lambda)=\lambda( R_1+J) + (R_0+a J)$, where $a\geq0$. Then the eigenvalues of $P(\lambda)$ are contained in the closed left half-plane. 
\end{corollary}

\begin{proof}
The case $a=0$ was considered in \cite{MehMW22}. Now take $a>0$ and define a matrix polynomial $P(\lambda)=\lambda^2\frac{R_1}a+\lambda \cdot 2 J+ R_0$, note that it clearly satisfies the assumptions of Theorem~\ref{half-plane}. Hence it is hyperstable, and by Theorem~\ref{Tkappa2} the matrix polynomial 
$$
(T_2 P)(z_1,z_2)=z_1z_2 \frac{R_1}a+(z_1+z_2) J+ R_0
$$
is hyperstable with respect to $H_{\pi/2}^2$. In particular,  if we set $z_2=a$ and replace $z_1$ by $\lambda$ then we obtain the original polynomial 
$\lambda a \frac{R_1}a+ (\lambda+a)J+ R_0=\lambda A_1+A_0$
being stable with respect to $H_{\pi/2}$.
\end{proof}

\section{Acknowledgment} Both authors acknowledge the financial support of the Priority Research Area SciMat
under the program Excellence Initiative Research University at
the Jagiellonian University in Krakow, Poland, decision no.  U1U/P05/NO/03.55.

\bibliography{szw}

\bibliographystyle{plain}

Faculty of Mathematics and Computer Science, Jagiellonian University, 
\L ojasiewicza 6, 30-348 Kraków, Poland  \\
E-mails : oskar.szymanski@doctoral.uj.edu.pl, michal.wojtylak@uj.edu.pl

\end{document}